\documentclass[11pt]{article}

\usepackage{amsmath,amsthm,amsfonts,amssymb,amscd}
\usepackage{subfigure}
\usepackage{graphicx,bm,color}
\usepackage{algorithm}
\usepackage{algpseudocode}
\usepackage{stmaryrd}
\usepackage{listings}
\usepackage{lineno}
\usepackage{mathrsfs}
\usepackage{lipsum}
\usepackage{url}            
\usepackage{booktabs}
\usepackage{authblk}
\usepackage[font=small,labelfont=bf]{caption}
\usepackage[justification=centering]{caption}

\usepackage[numbers]{natbib}

\definecolor{gray}{gray}{0.6}

\begin{document}

\title{Sequential Refinement Solver using Space-Time Domain Decomposition for Non-linear Multiphase Flow Problems}

\author{Hanyu Li and Mary F. Wheeler}
\date{\today}
\maketitle
\begin{abstract}
Convergence failure and slow convergence rate are among the biggest challenges with solving the system of non-linear equations numerically. While using strictly small time steps sizes and unconditionally stable fully implicit scheme mitigate the problem, the computational load becomes enormous. We introduce a sequential local refinement scheme in space-time domain that improves convergence rate and prevents convergence failure while not restricting to small time step, thus boosting computational efficiency. We rely on the non-linear two-phase flow model. The algorithm starts by solving the coarsest mesh. Then regions with certain features such as saturation front is refined to the finest resolution sequentially. Such process prevents convergence failure. After each refinement, the solution from the previous mesh is used to estimate initial guess of the current mesh for faster convergence. Numerical results are presented to confirm accuracy of our algorithm as compared to the traditional fine time step approach. We also observe 5 times speedup in the runtime by using our algorithm.
\end{abstract}
\newcommand{\bs}[1]{\boldsymbol{#1}}

\section{Introduction}
\label{sec:int}

\par Complex multi-phase flow and reactive transport in subsurface porous media is modeled by a system of non-linear equations. A common practice to solve such non-linear system is to approximate it in linear form and use iterative methods, such as Newton's method, to find the true solution. For large-scale models, such approach is usually computationally prohibitive even after parallelization. Due to the large number of unknowns, the approximate linear system becomes computationally exhaustive. More importantly, the significant non-linearity in the true system either requires a large number of iterations for convergence or results in failure of convergence when time-stepping is too aggressive. If the iterative method could be optimized such that, the number of iterations is minimized and the convergence is guaranteed, then we can achieve orders of magnitude greater computational efficiency.
\par Prior work exists to improve computational efficiency by reducing the size of the approximate linear system. Adaptive homogenization \cite{Amanbek:17, Singh:1118} addresses the problem by replacing fine grid with coarse grid in regions where non-linearity and variable (eg. saturation) variation is negligible, thus reducing the total number of spatial unknowns. However, fine and coarse grid in space requires different time scales for stable numerical solution. Forcing the coarse grid to accommodate the fine grid by taking fine time steps fails to reduce the number of unknowns in time. Space-time domain decomposition addresses this issue by allowing different time scales for different spatial grid, thus reducing the number of temporal unknowns. Several space-time domain decomposition approaches has been proposed in the past. \cite{Hughes:0288, Hulbert:1290} proposed space-time finite element method for elastodynamics with discontinuous Galerkin (DG) in time. The method has also been applied to other types of problems such as diffusion with different time discretization schemes \cite{Bause:1215, Bause:0617, Kocher:15, Kocher:1114}.
\par
\begin{figure}[H]
  \center{\includegraphics[width=\textwidth]{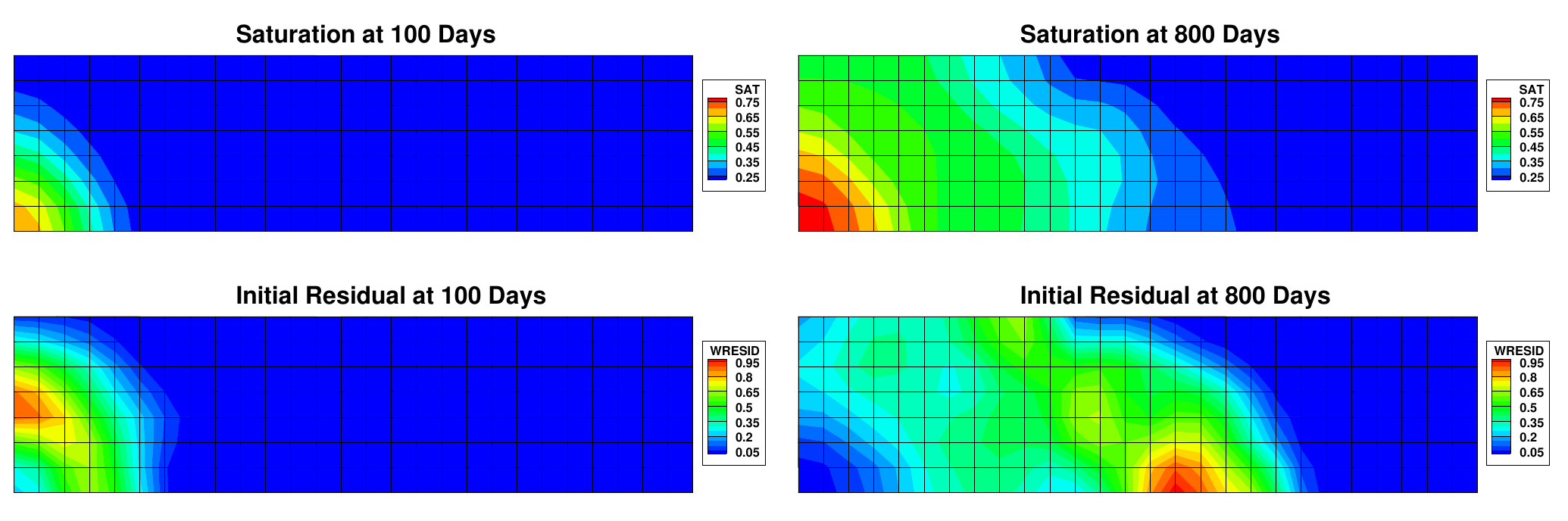}}
  \setlength{\abovecaptionskip}{-10pt}
  \caption{\bf{Saturation and normalized initial non-linear residual at 100 and 800 days}}
  \label{fig:satres}
\end{figure}

 The aforementioned literatures applied space-time decomposition method to mechanics problems. On the other hand, prior work regarding flow mostly focused on linear single phase flow and transport problems where flow is naturally decoupled from the advection-diffusion component transport \cite{Hoang:1213, Hoang:0717}. \cite{Singh:0818} first presented results for solving non-linear coupled multiphase flow and transport problem using space-time domain decomposition. \cite{Singh:0818} enforces strong continuity of fluxes at non-matching space-time interface with enhanced velocity. It also constructs and solves a monolithic system to avoid computational overheads associated with iterative solution schemes (\cite{Hoang:1213}) that require subdomain to be solved iteratively until weak continuity of fluxes is satisfied at interface. \cite{Singh:0918} further improves the method by allowing adaptive mesh refinement, thus improving computational efficiency while maintaining accuracy as compared to fine scale solution. It uses initial residual to search for regions that need refinement in space-time domain. As shown in Fig.\ref{fig:satres}, the normalized non-linear residual becomes the largest in the region with the highest non-linearity (saturation front) and thus consumes most computational resources and affects accuracy the most. Refining such region will reduce computational cost while maintaining accuracy as compared to solving the fine scale problem.
\par The adaptive local mesh refinement improves computational efficiency by reducing the size of the approximate linear system. \cite{Singh:0918} demonstrated the approach with only one level of refinement in both space and time, restricting the largest coarse time step allowed for stable numerical convergence. Also, no effort has been made to optimize the iterative method by reducing the number of iterations required for convergence. The iterative method approaches the true solution from the initial guess in a stepwise fashion. The rate of convergence in Newton's method heavily relies on the non-linearity possessed by the model and the initial guess. If the derivative of the residual function changes direction rapidly during iterations, it’s most likely to cause convergence failure. Meanwhile, if the initial guess is already close to the true solution, not many iterations are required to achieve convergence. Since derivative tends to be stable near the true solution, having a close initial guess becomes the key to avoid convergence failure and to improve convergence rate. In this work, we present a framework that allows several levels of refinement in space-time domain to represent features (eg. saturation front) of the system with the minimum number of grid cells. We will also optimize the convergence rate of the iterative method by providing better initial guess through sequential refinement.
\par In this work, we restrict ourselves to non-linear two-phase flow problems in subsurface porous media. We intent to approach more complicated non-linear problems such as black oil problem in the near future. The rest of the paper begins by describing the governing equations for two phase flow and its fully discrete form in Section \ref{sec:two}. Then we will present the solution algorithm for the sequential solver in Section \ref{sec:sol}. Afterwards, we demonstrate results from numerical experiments using the proposed algorithm in Section \ref{sec:num}.

\section{Two phase flow formulation}
\label{sec:two}

\subsection{Governing equations}
\label{subsec:gov}

\par We consider the following well-known two-phase, slightly compressible flow in porous medium model, with oil and water phase mass conservation, constitutive equations, boundary and initial conditions.
\begin{equation}
  \frac{\partial (\phi\rho_{\alpha}s_{\alpha})}{\partial t}+\nabla\cdot\bs{u}_{\alpha}=q_{\alpha}\quad in\ \Omega\, \times\, J
  \label{eq:mb}
\end{equation}
\begin{equation}
  \bs{u}_{\alpha}=-K\rho_{\alpha}\frac{k_{r\alpha}}{\mu_{\alpha}}(\nabla p_{\alpha}-\rho_{\alpha}\bs{g})\quad in\ \Omega\, \times\, J
  \label{eq:da}
\end{equation}
\begin{equation}
  \bs{u}_{\alpha}\cdot\bs{\nu}=0\quad on\ \partial\Omega\, \times\, J
  \label{eq:bc}
\end{equation}
\begin{equation}
  \begin{cases}
    p_{\alpha}=p_{\alpha}^{0} \\
    s_{\alpha}=s_{\alpha}^{0}
  \end{cases}
  at\ \Omega\, \times\, \{t=0\}
  \label{eq:ic}
\end{equation}
$\phi$ and $K$ are porosity and permeability tensor. $\rho_{\alpha}$, $s_{\alpha}$, $\bs{u}_{\alpha}$ and $q_{\alpha}$ are density, saturation, velocity and source/sink, respectively for each phase. The phases are slightly compressible and the phase densities are calculated by \eqref{eq:comp}.
\begin{equation}
  \rho_{\alpha}=\rho_{\alpha,ref}\cdot e^{c_{f,\alpha}(p_{\alpha}-p_{\alpha,ref})}
  \label{eq:comp}
\end{equation}
with $c_{f,\alpha}$ being the fluid compressibility and $\rho_{\alpha,ref}$ being the reference density at reference pressure $p_{\alpha,ref}$. In the constitutive equation \eqref{eq:da} given by Darcy's law, $k_{r\alpha}$, $\mu_{\alpha}$ and $p_{\alpha}$ are the relative permeability, viscosity and pressure for each phase. Relative permeability is a function of saturation. Pressure differs between wetting phase and non-wetting phase because of capillary pressure which is also a function of saturation.
\begin{equation}
  k_{r\alpha}=f(s_{\alpha})
  \label{eq:relperm}
\end{equation}
\begin{equation}
  p_{c}=g(s_{\alpha})=p_{nw}-p_{w}
  \label{eq:cap}
\end{equation}
The saturation of all phases obeys the constrain \eqref{eq:satcon}.
\begin{equation}
  \sum_{\alpha} s_{\alpha}=1
  \label{eq:satcon}
\end{equation}
The boundary and initial conditions are given by \eqref{eq:bc} and \eqref{eq:ic}. $J=(0,T]$ is the time domain of interest while $\Omega$ is the spatial domain.
\par Now we will give a brief introduction of mixed weak formulation in space-time domain. The functional spaces for mixed weak formulation are \\
\indent $\bs{V}=H(div;\Omega)=\Big\{\bs{v}\in\big(L^2(\Omega)\big)^d:\nabla\cdot\bs{v}\in L^2(\Omega)\Big\}$, \\
\indent $W=L^2(\Omega)$, \\
with finite dimensional subspace as $\bs{V}_{h}$ and $W_{h}$. As described in \cite{Singh:0818}, following the discontinuous Galerkin discretization in time, define space \\
\indent $\bs{V}_{h}^{t}=\bigg\{\bs{v}:J\rightarrow\bs{V}_{h}:\bs{v}\big|_{J_{m}}\in \big(\bs{P}_{l}(J_{m}) \big)^{d},m=1,\ldots,q,d=1,2\ or\ 3\bigg\}$, \\
\indent $W_{h}^{t}=\Big\{w:J\rightarrow W_{h}:w\big|_{J_{m}}\in \bs{P}_{l}(J_{m}),m=1,\ldots,q\Big\}$, \\
\indent $\bs{P}_{l}(J_{m})=\Bigg\{w:J_{m}\rightarrow W_{h}:w(t)=\displaystyle\sum_{a=1}^{l} w_{a}t^{a}\ with\ w_{a}\in W_{h}\Bigg\}$. \\
$d$ is the dimension of spatial domain. $\bs{V}_{h}^{t}$ and $W_{h}^{t}$ are spaces of functions that map from time domain $J$ to $\bs{V}_{h}$ and $W_{h}$ for each time interval $J_{m}$. These functions are represented by polynomials with degree up to $l$. In our framework, we will use $DG_{0}$ (polynomial of degree zero) discretization in time. Then the space-time mixed finite element space is \\
\indent $\bs{V}_{h}^{t,*}=\bs{V}_{h}^{t} \bigcap H(div;\Omega)\times J$. \\
Consider the oil-water system, the expanded variational form of Eqn.\eqref{eq:mb} through \eqref{eq:ic} is: find $\bs{u}_{\alpha,h}^{t}\in \bs{V}_{h}^{t,*}$, $\tilde{\bs{u}}_{\alpha,h}^{t}\in \bs{V}_{h}^{t,*}$, $S_{w,h}^{t}\in W_{h}^{t}$, $p_{o,h}^{t}\in W_{h}^{t}$ such that
\begin{equation}
  \bigg(\frac{\partial}{\partial t}\phi\Big(\rho_{w}s_{w,h}^{t}+\rho_{o}(1-s_{w,h}^{t})\Big),w\bigg)+\bigg(\nabla\cdot\Big(\bs{u}_{w,h}^{t}+\bs{u}_{o,h}^{t}\Big),w\bigg)=\bigg(q_{w}+q_{o},w\bigg)
  \label{eq:tot}
\end{equation}
\begin{equation}
  \bigg(\frac{\partial}{\partial t}\Big(\phi\rho_{w}s_{w,h}^{t}\Big),w\bigg)+\bigg(\nabla\cdot\bs{u}_{w,h}^{t},w\bigg)=\bigg(q_{w},w\bigg)
  \label{eq:wat}
\end{equation}
\begin{equation}
  \Big(K^{-1}\tilde{\bs{u}}_{o,h}^{t},\bs{v}\Big)-\Big(p_{o,h}^{t},\nabla\cdot\bs{v}\Big)=0
  \label{eq:omfe}
\end{equation}
\begin{equation}
  \Big(K^{-1}\tilde{\bs{u}}_{w,h}^{t},\bs{v}\Big)-\Big(p_{w,h}^{t},\nabla\cdot\bs{v}\Big)=-\Big(p_{c},\nabla\cdot\bs{v}\Big)
  \label{eq:wmfe}
\end{equation}
\begin{equation}
  \big(\bs{u}_{\alpha,h}^{t},\bs{v}\big)=\big(\lambda_{\alpha}\tilde{\bs{u}}_{\alpha,h}^{t},\bs{v}\big)
  \label{eq:aux}
\end{equation}
with $w\in W$ and $\bs{v}\in\bs{V}$. The mobility ratio in \eqref{eq:aux} is defined as
\begin{equation}
  \lambda_{\alpha}=\frac{k_{r\alpha}\rho_{\alpha}}{\mu_{\alpha}}
  \label{eq:mob}
\end{equation}
The additional auxiliary phase fluxes $\tilde{\bs{u}}_{\alpha}$ is used to avoid inverting zero phase relative permeability \cite{Peszy:0306}. The oil saturation and water pressure are eliminated by the saturation constrain and the capillary pressure relation (assume oil phase being the non-wetting phase).

\subsection{Fully discrete formulation}
\label{subsec:ful}

\par We will start by stating the basis functions in $RT_{0}\times DG_{0}$ discretization scheme. The pressure and saturation are piecewise constants while velocity is piecewise linear.
\begin{equation}
  w_{i}^{m}=
  \begin{cases}
    1\quad on\ E_{i}^{m}=x_{i-\frac{1}{2}}\leq x\leq x_{i+\frac{1}{2}}\bigcap t^{m}<t\leq t^{m+1} \\
    0\quad otherwise
  \end{cases}
  \label{eq:cbase}
\end{equation}
\begin{equation}
  \bs{\varphi}_{i+\frac{1}{2}}^{m}=
  \begin{cases}
    \frac{x-x_{i-\frac{1}{2}}}{\big|E_{i}^{m}\big|}\quad on\ E_{i}^{m} \\
    \frac{x_{i+\frac{3}{2}}-x}{\big|E_{i+1}^{m}\big|}\quad on\ E_{i+1}^{m}
  \end{cases}
  \label{eq:lbase}
\end{equation}
The solution to Eqn.\eqref{eq:tot} through \eqref{eq:aux} can be written in discrete form using the basis functions as
\begin{equation}
  \begin{cases}
    p_{o}=\displaystyle\sum_{m=1}^{q}\sum_{i=1}^{r}P_{i}^{m}w_{i}^{m} \\
    s_{w}=\displaystyle\sum_{m=1}^{q}\sum_{i=1}^{r}S_{w,i}^{m}w_{i}^{m} \\
    \bs{u}_{\alpha}=\displaystyle\sum_{m=1}^{q}\sum_{i=1}^{r+1}U_{\alpha,i+\frac{1}{2}}^{m}\bs{\varphi}_{i+\frac{1}{2}}^{m} \\
    \tilde{\bs{u}}_{\alpha}=\displaystyle\sum_{m=1}^{q}\sum_{i=1}^{r+1}\tilde{U}_{\alpha,i+\frac{1}{2}}^{m}\bs{\varphi}_{i+\frac{1}{2}}^{m} \\
  \end{cases}
  \label{eq:disol}
\end{equation}
We now substitute the testing functions in the variational forms of mass conservation and constitutive equation with $w_{j}^{n}$ and $\bs{\varphi}_{j+\frac{1}{2}}^{n}$, while keeping the solution in discrete form. For the first term in Eqn.\eqref{eq:omfe} and \eqref{eq:wmfe} we obtain
\begin{equation}
\begin{split}
\begin{aligned}
  \Big(K^{-1}\tilde{\bs{u}}_{\alpha},\bs{\varphi}_{j+\frac{1}{2}}^{n}\Big)_{\Omega\times J}&=\Bigg(K^{-1}\displaystyle\sum_{m=1}^{q}\sum_{i=1}^{r+1}\tilde{U}_{\alpha,i+\frac{1}{2}}^{m}\bs{\varphi}_{i+\frac{1}{2}}^{m},\bs{\varphi}_{j+\frac{1}{2}}^{n}\Bigg)_{\Omega\times J} \\ 
  &=\frac{1}{2\Big|e_{j+\frac{1}{2}}^n\Big|}\Bigg(\frac{x_{j+\frac{1}{2}}-x_{j-\frac{1}{2}}}{K_{j}}+\frac{x_{j+\frac{3}{2}}-x_{j+\frac{1}{2}}}{K_{j+1}}\Bigg)U_{\alpha,j+\frac{1}{2}}^{n}
\end{aligned}
\end{split}
\label{eq:fir}
\end{equation}
Here, $\Big|e_{j+\frac{1}{2}}^n\Big|$ is an edge of a space-time element. Since the framework uses backward Euler scheme in time to avoid Courant-Fredricks-Levy condition, we have the construction
\begin{equation}
  \bs{\varphi}_{i+\frac{1}{2}}^{m}(e_{j+\frac{1}{2}}^{n})=
  \begin{cases}
    \frac{1}{\big|e_{j+\frac{1}{2}}^n\big|}& \quad as\ i=j\ and\ m=n \\
    0                                                      & \quad otherwise
  \end{cases}
  \label{eq:stcon}
\end{equation}
The second term in Eqn.\eqref{eq:omfe} and \eqref{eq:wmfe} can be written as
\begin{equation}
  \begin{split}
  \begin{aligned}
  \Big(p_{\alpha},\nabla\cdot\bs{\varphi}_{j+\frac{1}{2}}^{n}\Big)_{\Omega\times J}&=\Bigg(\displaystyle\sum_{m=1}^{q}\sum_{i=1}^{r}P_{\alpha,i}^{m}w_{i}^{m},\nabla\cdot\bs{\varphi}_{j+\frac{1}{2}}^{n}\Bigg)_{\Omega\times J} \\
  &=\int\limits_{E_{j}^{n}}\frac{P_{\alpha,j}^{n}}{\big|E_{j}^{n}\big|}-\int\limits_{E_{j+1}^{n}}\frac{P_{\alpha,j+1}^{n}}{\big|E_{j+1}^{n}\big|} \\
  &=P_{\alpha,j}^{n}-P_{\alpha,j+1}^{n}
  \end{aligned}
  \end{split}
  \label{eq:sec}
\end{equation}
In case non-matching grid is encountered when the time scale is different at $(j+\frac{1}{2})^-$ and $(j+\frac{1}{2})^+$, assume the ratio between coarse and fine time step is $\frac{\delta t_{c}}{\delta t_{f}}=\tau$, then
\begin{equation}
  \begin{split}
  \begin{aligned}
    \Big(p_{\alpha},\nabla\cdot\bs{\varphi}_{j+\frac{1}{2}}^{n-\frac{1}{\tau}k}\Big)_{\Omega\times J}&=\Bigg(\displaystyle\sum_{m=1}^{q}\sum_{i=1}^{r}P_{\alpha,i}^{m}w_{i}^{m},\nabla\cdot\bs{\varphi}_{j+\frac{1}{2}}^{n-\frac{1}{\tau}k}\Bigg)_{\Omega\times J} \\
    &=P_{\alpha,j}^{n-\frac{1}{\tau}k}-P_{\alpha,j+1}^{n}
  \end{aligned}
  \end{split}
  \label{eq:nm}
\end{equation}
The variational form of capillary pressure term can be re-written in similar way as Eqn.\eqref{eq:sec} and \eqref{eq:nm}. Now we evaluate the mass conservation equation. The first term in Eqn.\eqref{eq:wat} becomes
\begin{equation}
  \bigg(\frac{\partial}{\partial t}\displaystyle\sum_{m=1}^{q}\sum_{i=1}^{r}\phi\rho_{w}s_{w,i}^{m}w_{i}^{m},w_{j}^{n}\bigg)_{\Omega\times J}=\Big((\phi\rho_{w}S_{w})_{j}^{n}-(\phi\rho_{w}S_{w})_{j}^{n-1}\Big)\big|E_{j}^{n-1}\big|
  \label{eq:wfir}
\end{equation}
In fine time scales, Eqn.\eqref{eq:wfir} can be altered as follow.
\begin{equation}
  \bigg(\frac{\partial}{\partial t}\displaystyle\sum_{m=1}^{q}\sum_{i=1}^{r}\phi\rho_{w}s_{w,i}^{m}w_{i}^{m},w_{j}^{n-\frac{1}{\tau}k}\bigg)_{\Omega\times J}=\Big((\phi\rho_{w}S_{w})_{j}^{n-\frac{1}{\tau}k}-(\phi\rho_{w}S_{w})_{j}^{n-\frac{1}{\tau}(k+1)}\Big)\Big|E_{j}^{n-\frac{1}{\tau}(k+1)}\Big|
  \label{eq:wfin}
\end{equation}
The second term is calculated as
\begin{equation}
  \begin{split}
  \begin{aligned}
    (\nabla\cdot\bs{u}_{w},w_{j}^{n})_{\Omega\times J}&=\bigg(\nabla\cdot\displaystyle\sum_{m=1}^{q}\sum_{i=1}^{r+1}U_{\alpha,i+\frac{1}{2}}^{m}\bs{\varphi}_{i+\frac{1}{2}}^{m},w_{j}^{n}\bigg)_{\Omega\times J} \\
    &=U_{w,j+\frac{1}{2}}^{n}-U_{w,j-\frac{1}{2}}^{n}
  \end{aligned}
  \end{split}
   \label{eq:wsec}
\end{equation}
The situation for non-matching grid is a little different for this term. Assume fine time stays on $(j+\frac{1}{2})^-$ side, then on the fine time element we have
\begin{equation}
  (\nabla\cdot\bs{u}_w,w_{j}^{n-\frac{1}{\tau}k})=U_{w,j+\frac{1}{2}}^{n-\frac{1}{\tau}k}-U_{w,j-\frac{1}{2}}^{n-\frac{1}{\tau}k}
  \label{eq:wsecf}
\end{equation}
while for coarse time element we have
\begin{equation}
  (\nabla\cdot\bs{u}_w,w_{j+1}^{n})=U_{w,j+\frac{3}{2}}^{n}-\displaystyle\sum_{k=1}^{\tau}U_{w,j+\frac{1}{2}}^{n-\frac{1}{\tau}k}
  \label{eq:wsecc}
\end{equation}
The oil phase mass conservation equation is similar. Adding the equation for these two phases will provide the expression for the total mass conservation equation. The two sides of Eqn.\eqref{eq:aux} is estimated as
\begin{equation}
  \begin{split}
  \begin{aligned}
    (\bs{u}_{\alpha},\bs{v})&=\displaystyle\sum_{m=1}^{q}\sum_{i=1}^{r+1}U_{\alpha,i+\frac{1}{2}}^{m}\Big(\bs{\varphi}_{i+\frac{1}{2}}^{m},\bs{\varphi}_{j+\frac{1}{2}}^{n}\Big) \\
    &=\frac{x_{j+\frac{3}{2}}-x_{j-\frac{1}{2}}}{2\big|e_{j+\frac{1}{2}}^{n}\big|}U_{\alpha,j+\frac{1}{2}}^{n}
  \end{aligned}
  \end{split}
  \label{eq:auxl}
\end{equation}
\begin{equation}
  (\lambda_{\alpha}\tilde{\bs{u}}_{\alpha},\bs{v})\approx(\lambda_{\alpha}^{*}\tilde{\bs{u}}_{\alpha},\bs{v})=\frac{x_{j+\frac{3}{2}}-x_{j-\frac{1}{2}}}{2\big|e_{j+\frac{1}{2}}^{n}\big|}\lambda_{\alpha,j+\frac{1}{2}}^{*,n}U_{\alpha,j+\frac{1}{2}}^{n}
  \label{eq:auxr}
\end{equation}
The $\lambda_{\alpha,j+\frac{1}{2}}^{*,n}$ is the upwind mobility for stable numerical solution and is defined as
\begin{equation}
  \lambda_{\alpha,j+\frac{1}{2}}^{*,n}=\rho_{\alpha,j+\frac{1}{2}}^{*,n}\frac{k_{j+\frac{1}{2}}^{r\alpha,*}}{\mu_{\alpha}}=
  \begin{cases}
    \frac{1}{2\mu_{\alpha}}(\rho_{\alpha,j}^{n}+\rho_{\alpha,j+1}^{n})k_{r\alpha}(S_{\alpha,j}^{n})     &\quad if\ \tilde{U}_{\alpha,j+\frac{1}{2}}^{n}>0 \\
    \frac{1}{2\mu_{\alpha}}(\rho_{\alpha,j}^{n}+\rho_{\alpha,j+1}^{n})k_{r\alpha}(S_{\alpha,j+1}^{n}) &\quad otherwise
  \end{cases}
  \label{eq:uw}
\end{equation}
The above section provides us a non-linear system of equations of pressure and saturation. To solve such system, we linearize it and use Newton's method to approach the true solution through iteration process. Depending on the level of non-linearity and the closeness between initial guess and true solution, Newton's method could take numerous iterations before achieving convergence. In the next section, we will introduce our sequential local refinement algorithm to minimize the number of iterations while maintaining solution accuracy.

\section{Solution algorithm}
\label{sec:sol}

\subsection{Sequential local refinement}
\label{subsec:seq}

\par
\begin{figure}
  \center{\includegraphics[width=\textwidth]{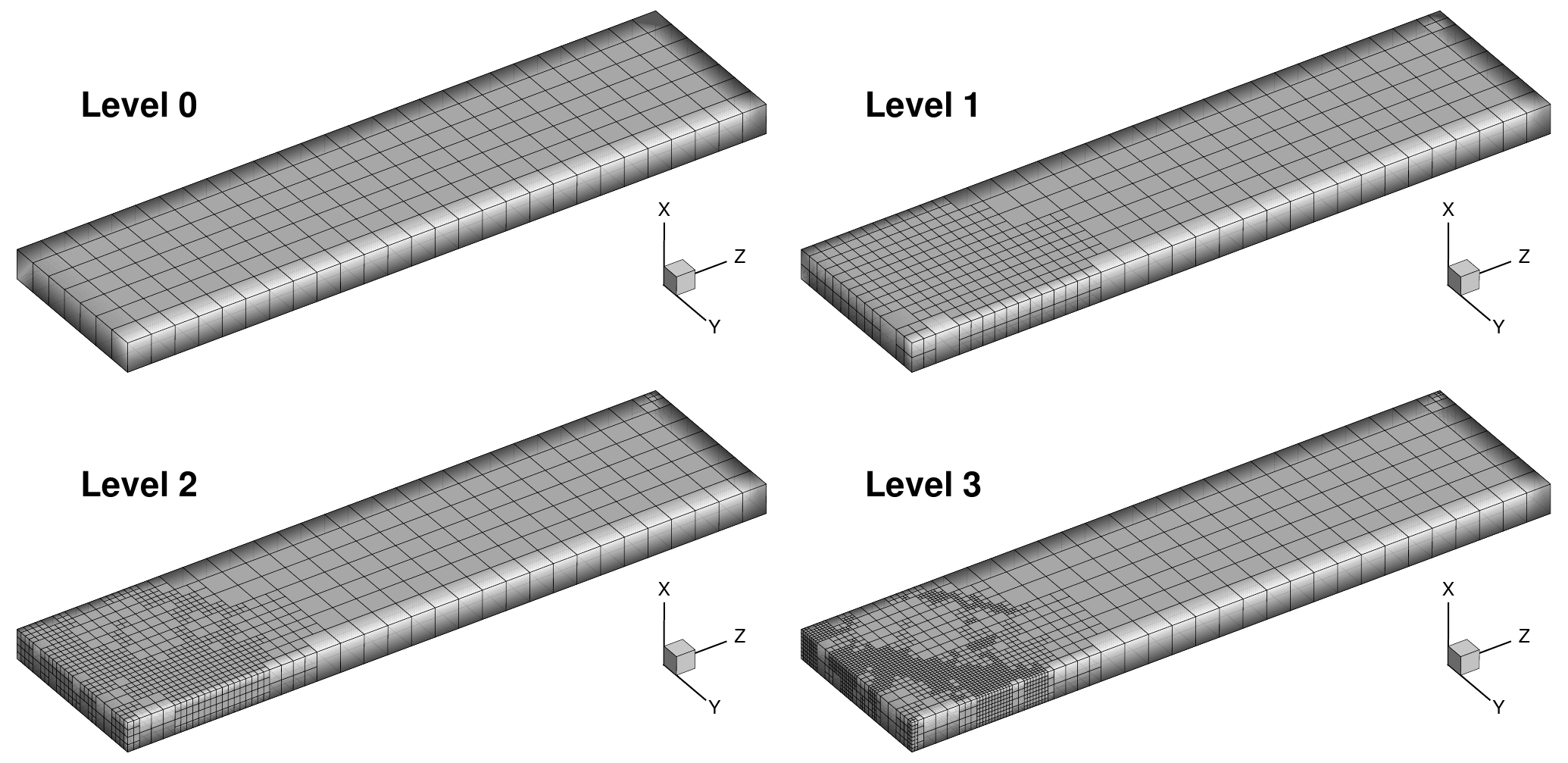}}
  \setlength{\abovecaptionskip}{-10pt}
  \caption{\bf{Sequential local mesh refinement in space-time domain from coarsest (level 0) to finest (level 1) resolution}}
  \label{fig:grid}
\end{figure}

 In this section we present the solver algorithm that uses sequential local refinement in space-time domain to provide initial guess close to the true solution, thus reducing the time for Newton convergence. The algorithm starts by solving the problem at its coarsest resolution. Then the given domain is sequentially refined isotropically in space-time domain to its finest resolution in regions colored by specific indicators. Fig.\ref{fig:grid} demonstrates a sample semi-structured grid generated during sequential local refinement. Here the $x$ axis represents time in 2-D spatial problem. Please note that we always refine cells that contain wells for accurate estimate of rate and bottom-hole pressure.
\par
\begin{figure}[H]
  \center{\includegraphics[width=\textwidth]{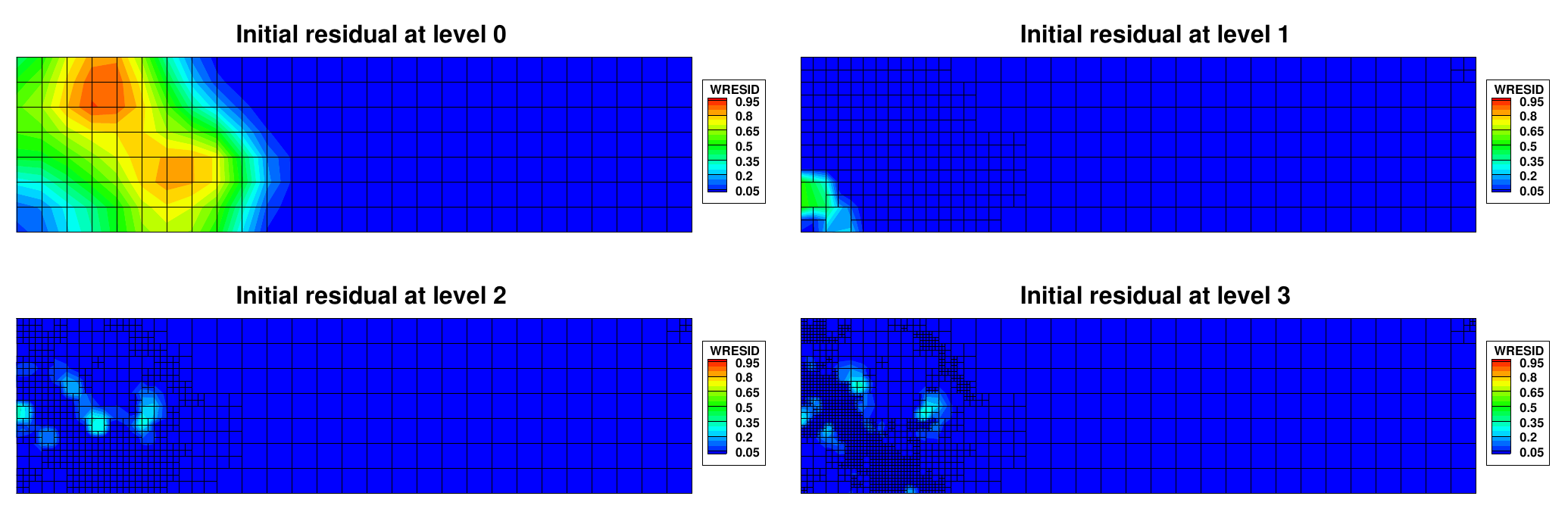}}
  \setlength{\abovecaptionskip}{-10pt}
  \caption{\bf{Normalized initial non-linear residual at each refinement level}}
  \label{fig:resid}
\end{figure}

\begin{figure}[H]
  \center{\includegraphics[width=\textwidth]{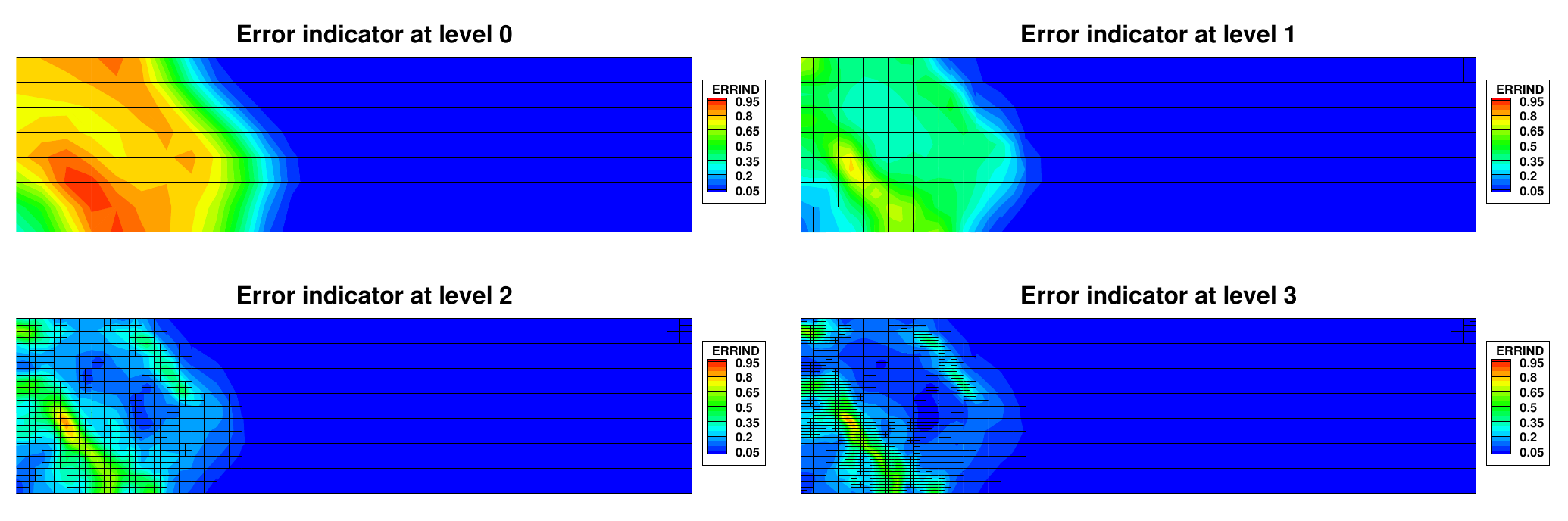}}
  \setlength{\abovecaptionskip}{-10pt}
  \caption{\bf{Normalized error indicator at each refinement level}}
  \label{fig:errind}
\end{figure}

\begin{figure}[H]
  \center{\includegraphics[width=\textwidth]{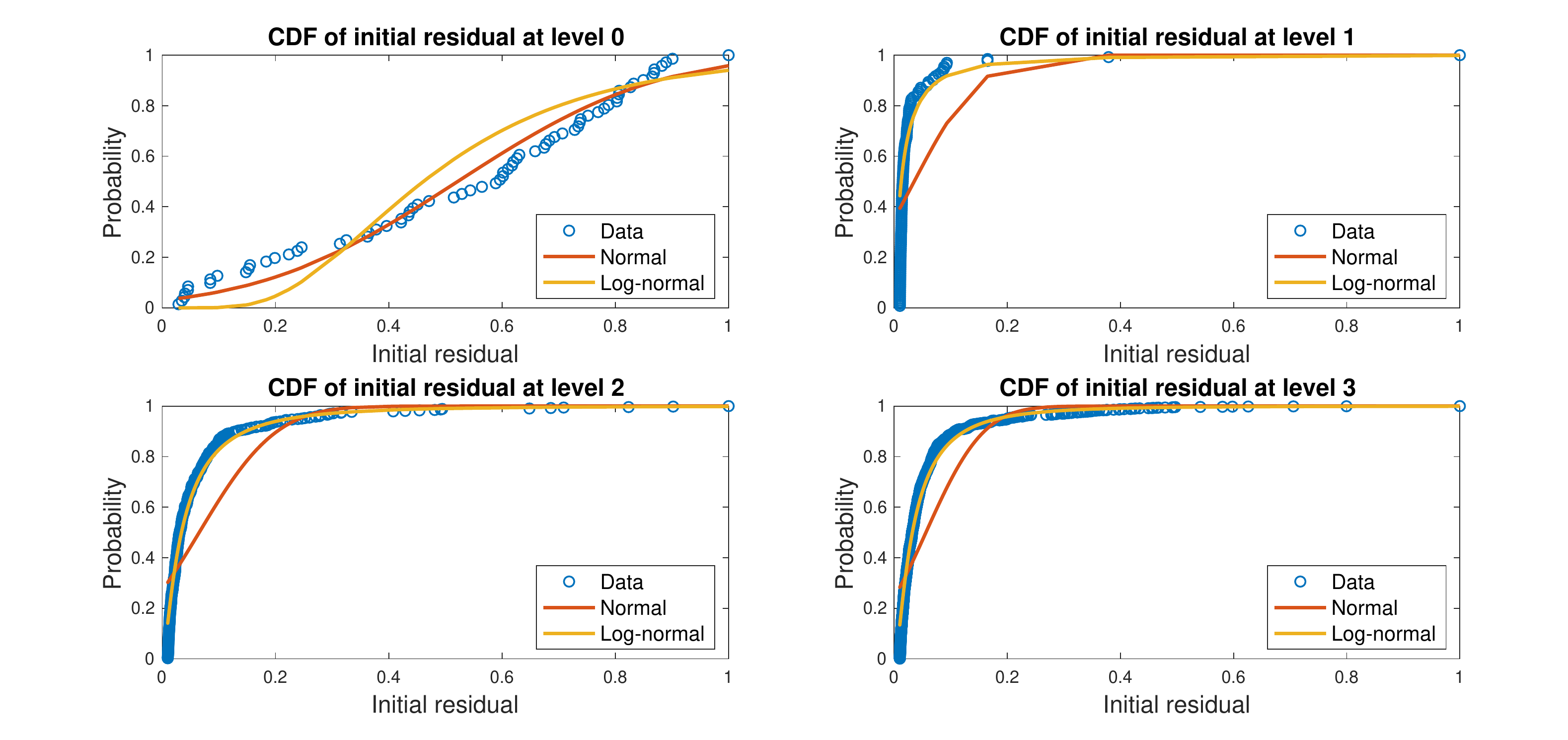}}
  \setlength{\abovecaptionskip}{-10pt}
  \caption{\bf{Cumulative distribution function fitted to initial residual data at each refinement level}}
  \label{fig:rcdf}
\end{figure}

\begin{figure}[H]
  \center{\includegraphics[width=\textwidth]{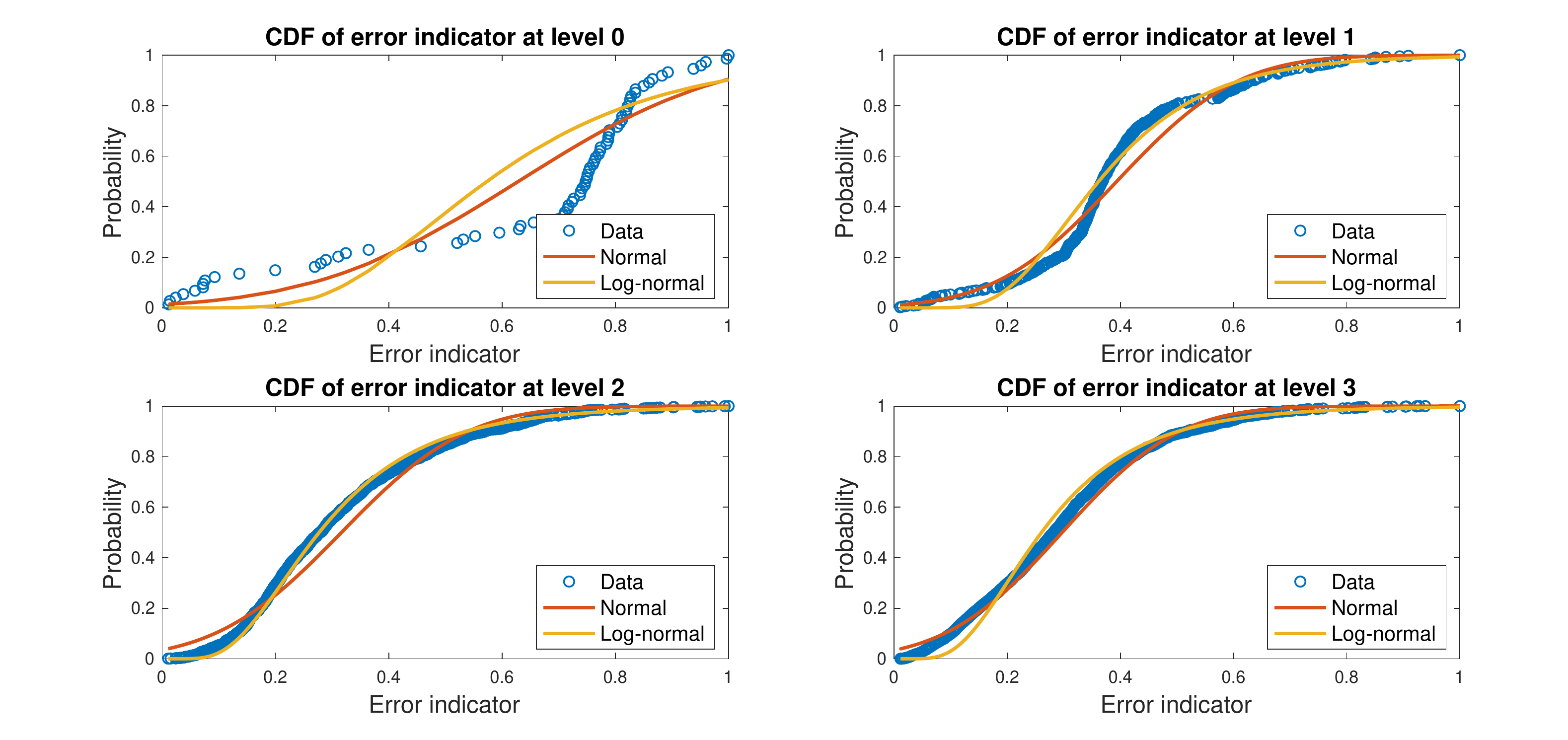}}
  \setlength{\abovecaptionskip}{-10pt}
  \caption{\bf{Cumulative distribution function fitted to error indicator data at each refinement level}}
  \label{fig:ecdf}
\end{figure}

\begin{figure}[H]
  \center{\includegraphics[scale=1]{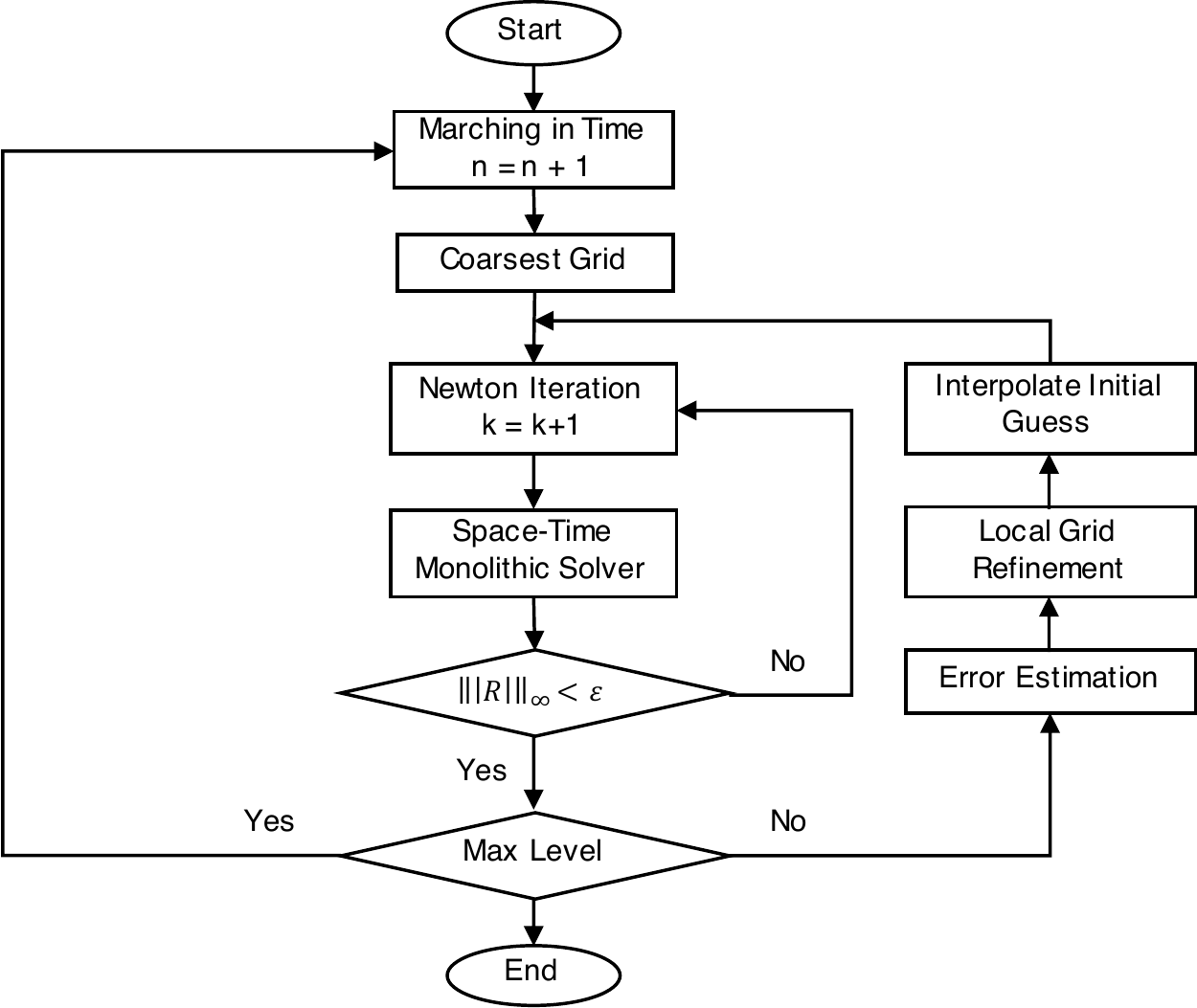}}
  \setlength{\abovecaptionskip}{5pt}
  \caption{\bf{Solution algorithm for sequential local refinement solver}}
  \label{fig:algo}
\end{figure}

 After each refinement, before solving the problem on the new mesh, the unknowns of newly generated fine elements are populated by the solution of the previous mesh using spatial and temporal linear interpolation. Such approach provides a close initial guess to the true solution, however it also creates a problem. The indicator used in \cite{Singh:0918} is the normalized initial non-linear residual calculated as
\begin{equation}
  \tilde{R}=\frac{|R|}{\big\||R|\big\|_{\infty}}
  \label{normre}
\end{equation}
It measures the closeness between the initial guess and the true solution. This indicator works perfectly during single level refinement. However, for multiple level refinement, since we are providing initial guess through linear interpolation, the initial guess on refined grid is naturally closer to the true solution and thus the initial residual does not expose certain feature of the system anymore. Fig.\ref{fig:resid} provides an example. After the first level of refinement, the interpolation calculates very close initial guess that causes the residual to appear only sporadicly.
The observation on initial residuals infers that we need another indicator to track features at each refinement level. The changes of certain solution variable in space and time measures the solution sensitivity to different scales. A large change in unit space/time indicates the existence of feature and that refinement provides more accurate solution. So we define an error indicator as
\begin{equation}
  \varepsilon=\bigg(\frac{(\Delta_{s}S_{w})^{2}}{\|(\Delta_{s}S_{w})^{2}\|_{\infty}}+\frac{(\Delta_{t}S_{w})^{2}}{\|(\Delta_{t}S_{w})^{2}\|_{\infty}}\bigg)^{\frac{1}{2}}
  \label{eq:errind}
\end{equation}
\begin{equation}
  \tilde{\varepsilon}=\frac{\varepsilon}{\| \varepsilon \|_{\infty}}
  \label{eq:normer}
\end{equation}
Here we use the saturation to calculate the error indicator because the pressure solution is too smooth. The change in space and time are normalized respectively so that they will have the same weight on calculating the error indicator. Eqn.\eqref{eq:normer} simply normalizes the error indicator to $[0,1]$ scale. Fig.\ref{fig:errind} shows the normalized error indicator distribution that exposes feature at each refinement level.
Now the question is how do we choose the region for refinement. We first define the range $[0.01,1]$ of normalized residual and error indicator as the analysis range. Anything below $0.01$ is neglected. At each level, we will refine the region with $50\%$ of the largest normalized values in the analysis range. The cumulative distribution function of initial residual and error indicator at each refinement level is plotted in Fig.\ref{fig:rcdf} and Fig.\ref{fig:ecdf} against sample data recorded during simulation. As demonstrated by the graphs, the initial residual data is better represented by log-normal distribution. Meanwhile the error indicator data follows the trend between normal and log-normal distribution. Therefore, the threshold for initial residual is the log-mean while for error indicator is the average of mean and log-mean.
\noindent Please note that the sporadic appearance of large initial residuals after the first refinement is also reflected in its cumulative distribution function. In Fig.\ref{fig:rcdf} from level 1 to 3, the initial residual samples stay concentrated towards $0$, unlike the error indicator samples that spread smoothly across the range. The complete algorithm is illustrated in Fig.\ref{fig:algo}.

\subsection{Data structure}
\label{subsec:dat}
\par The local mesh refinement process creates semi-structured grid which is stored in a tree formation. Each element on the coarsest grid is represented by a root node. All the other nodes in the tree are created during grid refinement. Each node is linked to its parent and children by pointers. This data structure facilitates the sequential refinement process as we can simply evolve the tree instead of creating every refined grid from scratch. After grid generation is complete, all the elements are indexed to the construct the monolithic system for the solver.
\par
\begin{figure}[H]
  \center{\includegraphics[scale=1]{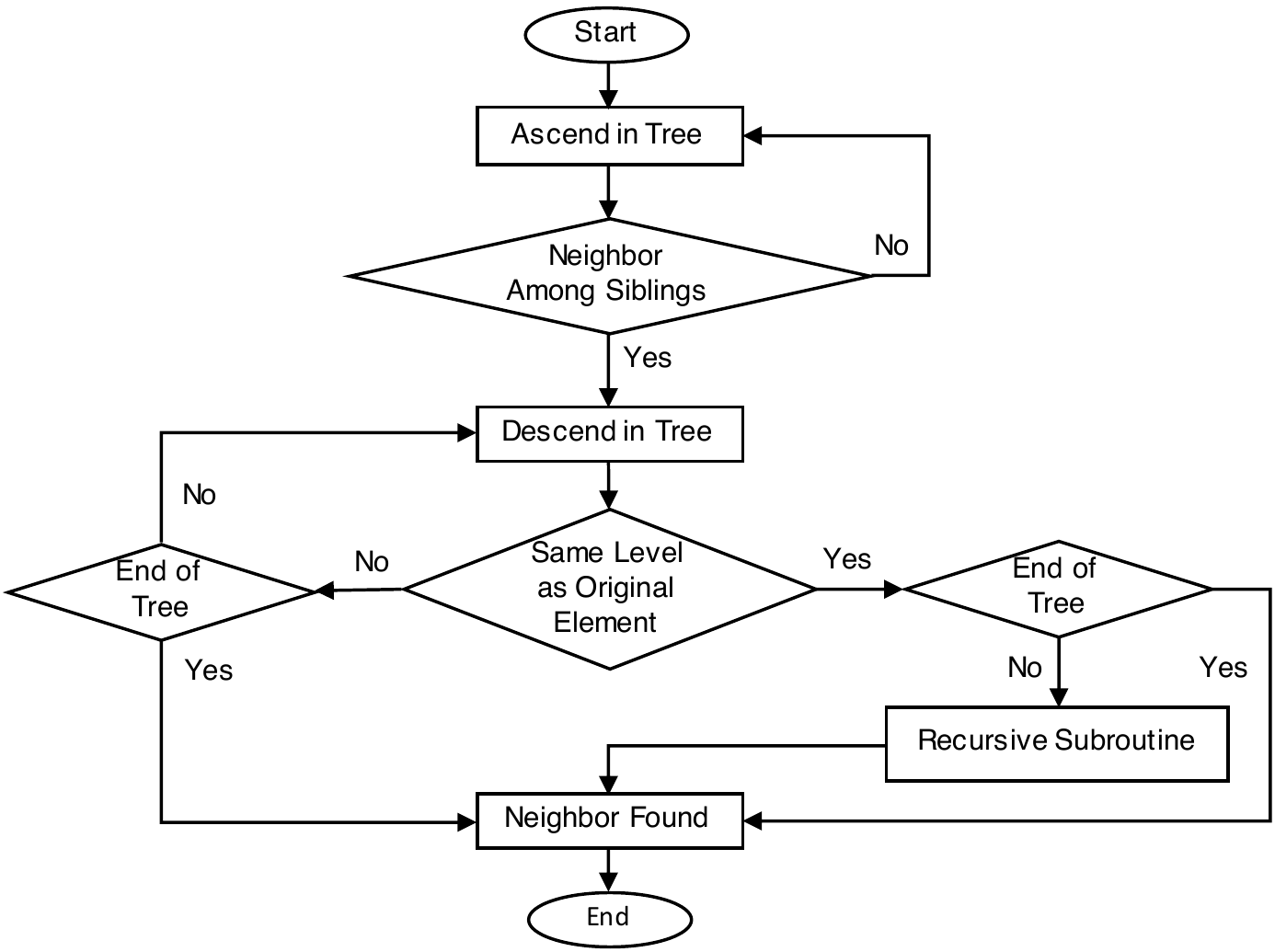}}
  \setlength{\abovecaptionskip}{5pt}
  \caption{\bf{Neighbor searching algorithm}}
  \label{fig:nalg}
\end{figure}

\begin{figure}[H]
  \center{\includegraphics[width=\textwidth]{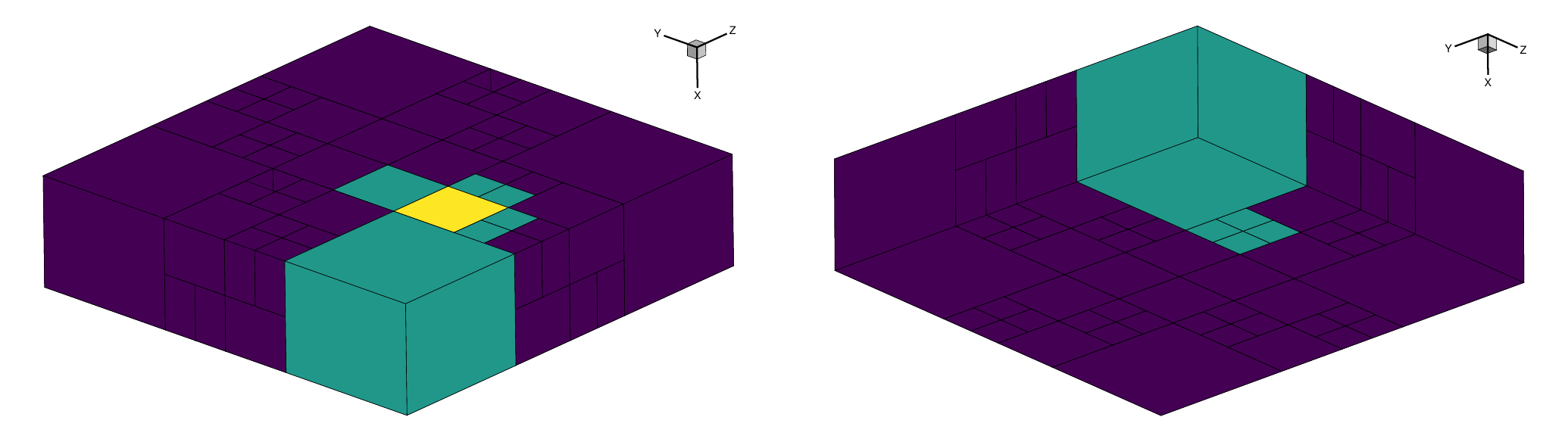}}
  \setlength{\abovecaptionskip}{-10pt}
  \caption{\bf{Neighbor searching algorithm outcome from two viewing angle}}
  \label{fig:neigh}
\end{figure}

 To successfully construct the monolithic system, we need to accurately pinpoint the neighbors given a specific element. We designed an algorithm to search the neighbors in each space-time direction separately (front/back/left/right/top/bottom/past/future. Future neighbor is not used for calculation. We only search it for auxiliary purposes such as visualization). During the search, we first ascend the tree from the original element until a neighbor exists among the sibling elements in the intended search direction. Then we move to that sibling and start descending the tree. If the descend terminates at or before the refinement level of the original element, then we have found the one and only neighbor. If more levels exist after descending to the same level as the original element, we use a recursive subroutine to locate all the neighbors on the deeper levels. The neighbor searching algorithm is shown in Fig.\ref{fig:nalg}.
\noindent We give an example of the neighbor searching algorithm outcome in Fig.\ref{fig:neigh} based on a complex semi-structured grid. Shown in two different angles, the original element is colored in yellow and the neighbor elements are colored in green.
\noindent Now that we have introduced the all the important algorithms involved in our solver, in the next section we will present some numerical results.

\section{Numerical results}
\label{sec:num}

\par
\begin{figure}[H]
  \center{\includegraphics[width=\textwidth]{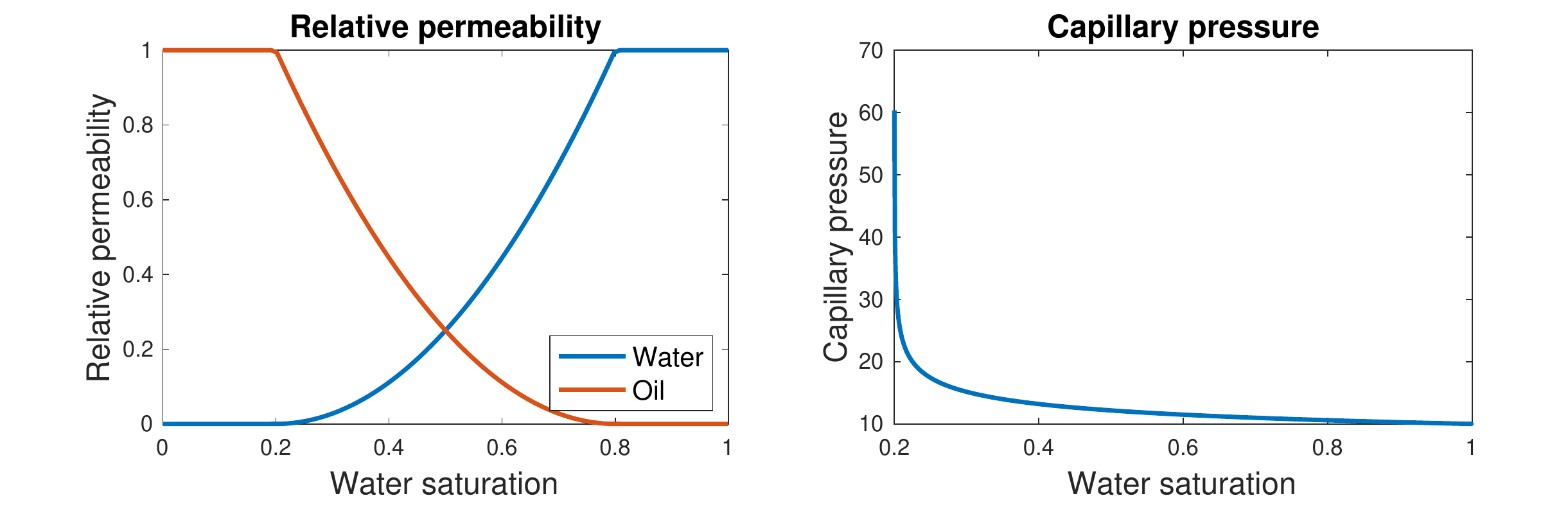}}
  \setlength{\abovecaptionskip}{-10pt}
  \caption{\bf{Relative permeability (left) and capillary pressure (right) curve for numerical experiments}}
  \label{fig:brooks}
\end{figure}

\begin{figure}[H]
  \center{\includegraphics[width=\textwidth]{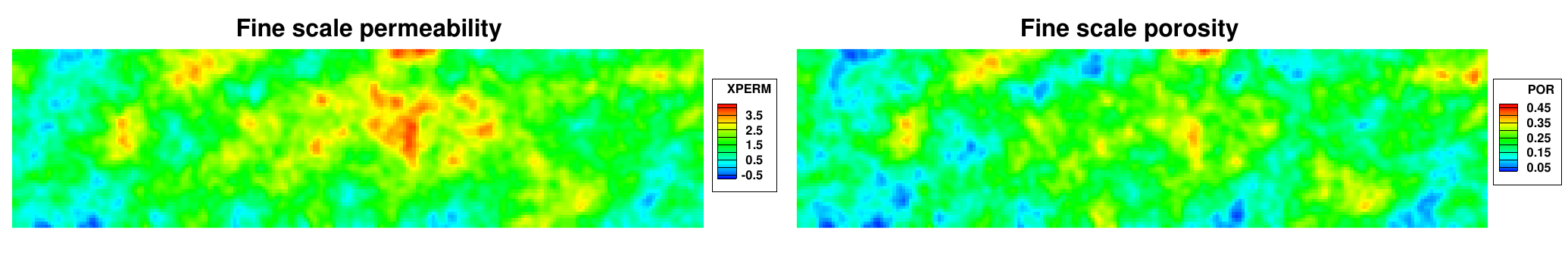}}
  \setlength{\abovecaptionskip}{-10pt}
  \caption{\bf{Gaussian-like fine permeability (left) and porosity (right) distribution}}
  \label{fig:gaup}
\end{figure}

\begin{figure}[H]
  \center{\includegraphics[width=\textwidth]{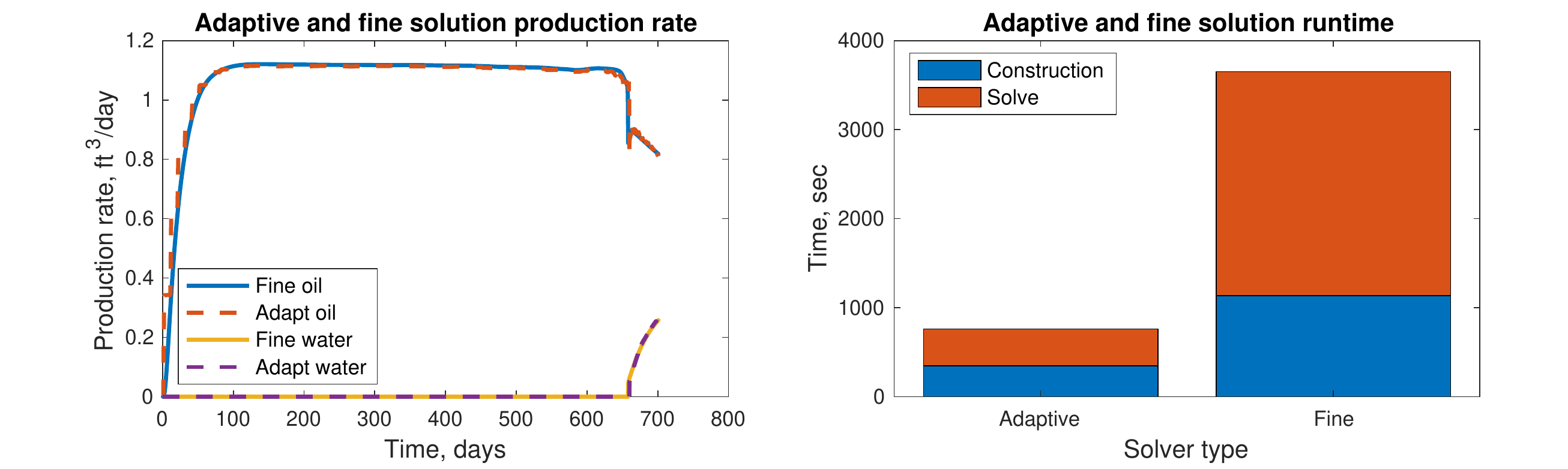}}
  \setlength{\abovecaptionskip}{-10pt}
  \caption{\bf{Two phase production rate and runtime from adaptive and fine solution}}
  \label{fig:gaura}
\end{figure}

 In this section we will show results from two numerical experiments with 2-D two phase flow problem. Both experiments use the same fluid data from the SPE10 (\cite{Christie:0801}) dataset. The oil and water reference densities in Eqn.\eqref{eq:comp} are taken to be $53\ [lb/ft^{3}]$ and $64\ [lb/ft^{3}]$ respectively and compressibilities are $1\times 10^{-4}\ psi^{-1}$ and $3\times 10^{-6}\ psi^{-1}$. We use Brooks's Corey model for both relative permeability and capillary pressure. The equations for relative permeability are
\begin{equation}
  \begin{cases}
    k_{rw}=k_{rw}^{0}\Big(\frac{S_{w}-S_{wirr}}{1-S_{or}-S_{wirr}}\Big)^{n_{w}} \\
    k_{ro}=k_{ro}^{0}\Big(\frac{S_{o}-S_{or}}{1-S_{or}-S_{wirr}}\Big)^{n_{o}}
  \end{cases}
  \label{eq:bcrel}
\end{equation}
The endpoint values are $S_{or}=S_{wirr}=0.2$ and $k_{ro}^{0}=k_{rw}^{0}=1.0$ and the model exponents are $n_{w}=n_{o}=2$. The equation for capillary pressure is
\begin{equation}
  p_{c}(S_{w})=P_{en,cow}\Big(\frac{1-S_{wirr}}{S_{w}-S_{wirr}}\Big)^{l_{cow}}
  \label{eq:bccap}
\end{equation}
with $P_{en,cow}=10\ psi$ and $l_{cow}=0.2$. Fig.\ref{fig:brooks} visualizes the relative permeability and capillary pressure curve.
In the first experiment, we use a gaussian-like permeability and porosity distribution which is plotted in Fig.\ref{fig:gaup}.
\noindent The computational domain is $56ft\times 216ft\times 1ft\times 700days$ with coarsest and finest element size of $8ft\times 8ft\times 1ft\times 10days$ and $1ft\times 1ft\times 1ft\times 1.25days$. We allow three refinement levels in our experiment. Although our framework allows each level to have different refinement ratio, for the sake of simplicity we will set the same ratio, a factor of 2, for all three levels. We place a rate specified injection well at bottom left corner and a pressure specified production well at upper right corner. The injection rate is $1\ ft^{3}/day$ and production pressure is $1000\ psi$. Furthermore, the initial pressure and water saturation are set to be $1000\ psi$ and $0.2$.
\begin{figure}[H]
  \center{\includegraphics[width=\textwidth]{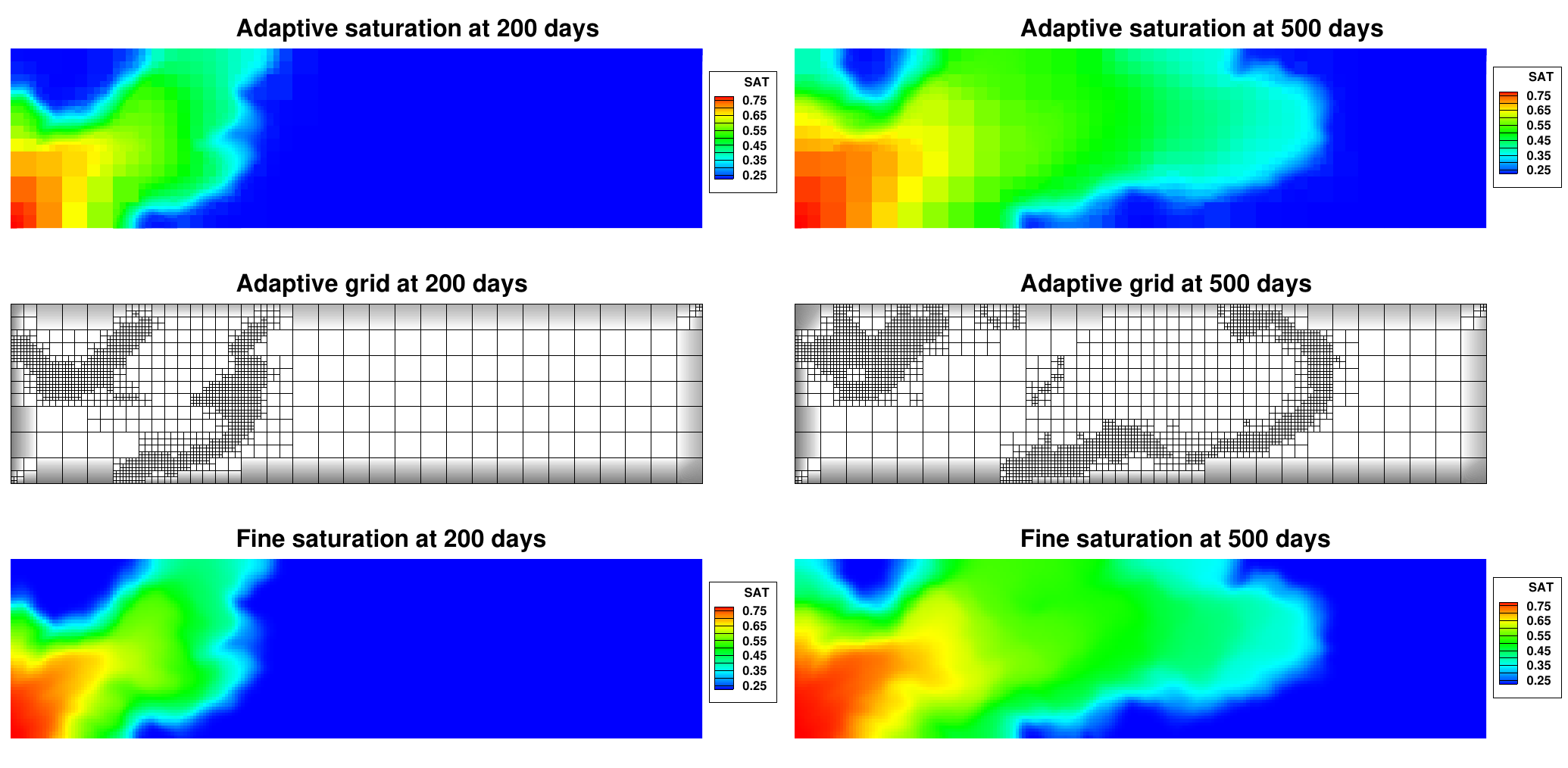}}
  \setlength{\abovecaptionskip}{-10pt}
  \caption{\bf{Adaptive saturation distribution (top) and adaptive mesh (middle) generated by sequential refinement solver as compared to fine scale solution (bottom) at 200 and 500 days}}
  \label{fig:gaure}
\end{figure}
Fig.\ref{fig:gaure} shows the adaptive grid saturation profile along with its mesh as compared to fine grid saturation profile at 200 and 500 days. The shape of the front looks similar while elements are coarsened behind the front in the adaptive solver. The production rates and runtime of the two solutions are plotted in Fig.\ref{fig:gaura}.
\noindent The oil and water production rate matches well between the two different solutions which indicates that the sequential solver is accurate. While maintaining the accuracy, the sequential solver using adaptive grid reduces the total runtime by approximately $5$ times with linear system construction and solving providing $4$ times and $6$ times reduction respectively.
\par Although the solver improves computational efficiency while maintaining accuracy, we observe several problems during the experiment. The first major problem is over refining the elements. The saturation has the biggest change in time at the front. Meanwhile, behind the front the saturation is stable time wise but varies significantly in space. Since we rely on isotropic refinement in space-time domain and our error indicator is calculated with both spatial and temporal variation of saturation, we are forcing the front region to be fine in space and behind-the-front region to be fine in time, which is redundant. Such over-refinement could cause sever increase in time needed for linear system construction. The isotropic refinement could also ignore some important features of the system. The second major problem is associated with estimating initial guess after each refinement. The isotropic refinement scheme requires us to interpolate spatially and temporally at the same time. The calculated initial guess does prevent convergence failure and improve convergence rate but not to a significant extent. These two problems suggest that separating time and space refinement is the solution.
\par
\begin{figure}[H]
  \center{\includegraphics[width=\textwidth]{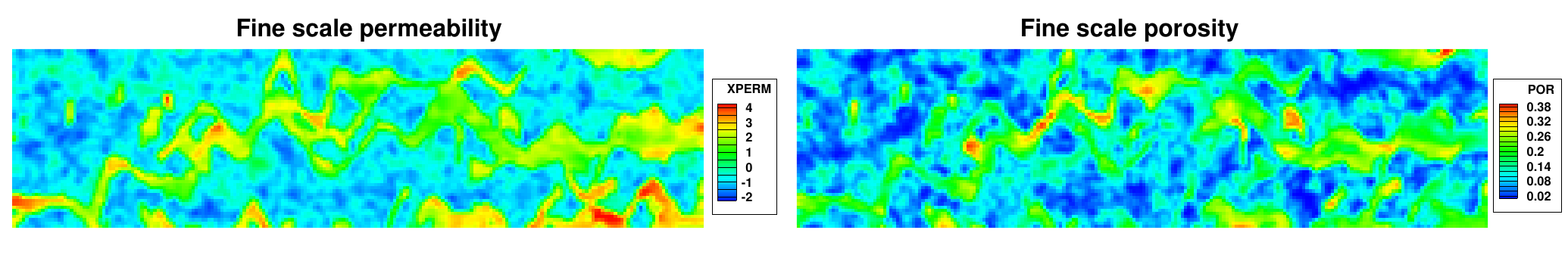}}
  \setlength{\abovecaptionskip}{-10pt}
  \caption{\bf{Channel-like fine permeability (left) and porosity (right) distribution}}
  \label{fig:chap}
\end{figure}

\begin{figure}[H]
  \center{\includegraphics[width=\textwidth]{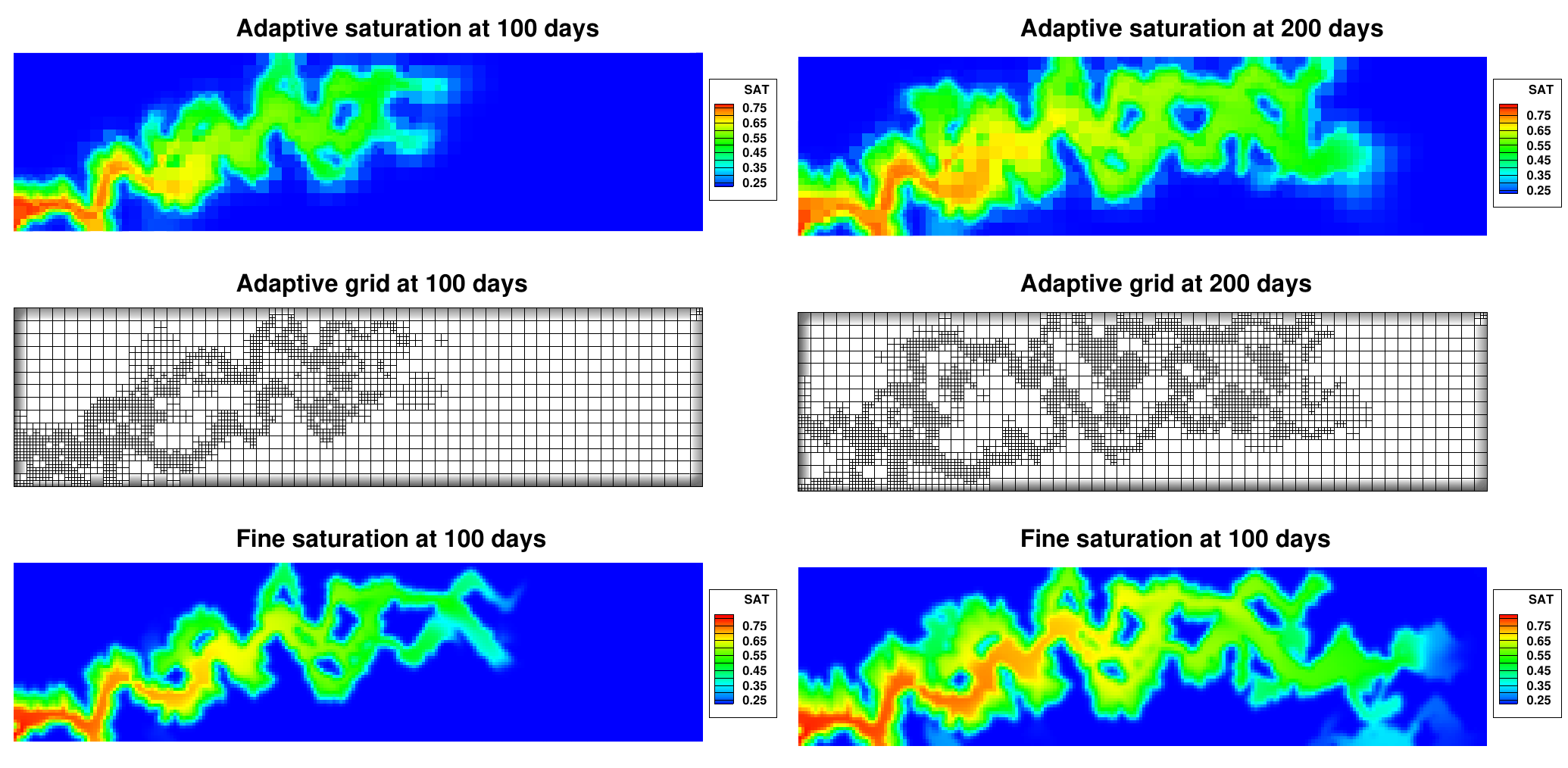}}
  \setlength{\abovecaptionskip}{-10pt}
  \caption{\bf{Adaptive saturation distribution (top) and adaptive mesh (middle) generated by sequential refinement solver as compared to fine scale solution (bottom) at 100 and 200 days}}
  \label{fig:chare}
\end{figure}

 In the second experiment, we use a channel-like permeability and porosity distribution as demonstrated by Fig.\ref{fig:chap}.
\noindent Considering the dramatic changes in petrophysical properties, we only allow two levels of refinement so that the coarsest level (level 0) does not destroy too much features of the system. For the same reason, the coarsest time scale has to be small to prevent convergence failure. The computational domain is $56ft\times 216ft\times 1ft\times 200days$ with coarsest and finest element size of $4ft\times 4ft\times 1ft\times 0.5day$ and $1ft\times 1ft\times 1ft\times 0.125day$. The injection rate, production pressure and initial condition is kept the same as the previous experiment. Fig.\ref{fig:chare} shows the adaptive grid saturation profile along with its mesh as compared to fine grid saturation profile at 100 and 200 days. The shape of the front looks similar however we do observe some minor feature lost in the coarse region of the adaptive grid. The algorithm also reduced the runtime by 5 times, however due to the lost of features, the production rates are less accurate. Meanwhile, the problems that appeared in the previous experiment deteriorates in this channel case. These observations also indicates that we need to separate the time and space refinement, which will be done next in the near future.

\section{Conclusions}
\label{sec:con}
\par We have presented an algorithm that sequentially refines the coarse mesh to solve non-linear two phase flow problems. After each refinement, the previous solution is used to interpolate the initial guess for the new mesh. Results from two numerical experiments are demonstrated. We have achieved 5 times speedup in computational time using our algorithm by both reducing the number elements and providing a better initial guess. Convergence failure is better prevented and convergence rate is improved. The result for the Gaussian-like permeability field shows similar saturation profile between the adaptive solution and fine solution. The production rates of the two solutions also match pretty well.  For the channel case the sequential solver becomes less accurate. The main problem that caused the inaccuracy is the isotropic refinement in space-time domain. The error indicator calculated by both spatial and temporal variation may sometimes mislead the refinement process and fail to capture certain features in the system. It also over-refine the coarse grid, causing increased runtime. We will solve this problem by separating the spatial and temporal refinement in the near future.

\bibliographystyle{plain} 
\bibliography{st_draft_1.bbl}

\begin{thebibliography}{10}

\bibitem{Amanbek:17}
Y.~Amanbek, G.~Singh, and M.F. Wheeler.
\newblock Adaptive numerical homogenization for upscaling single phase flow and
  transport.
\newblock {\em ICES Report}, 2017.

\bibitem{Bause:1215}
M.~Bause and U.~K{\"o}cher.
\newblock Variational time discretization for mixed finite element
  approximations of nonstationary diffusion problems.
\newblock {\em Journal of Computational and Applied Mathematics}, 289:208--224,
  December 2015.

\bibitem{Bause:0617}
M.~Bause, F.A. Radu, and U.~K{\"o}cher.
\newblock Space-time finite element approximation of the biot poroelasticity
  system with iterative coupling.
\newblock {\em Computer Methods in Applied Mechanics and Engineering},
  320:745--768, June 2017.

\bibitem{Christie:0801}
M.A. Christie and M.J. Blunt.
\newblock Tenth spe comparative solution project: A comparison of upscaling
  techniques.
\newblock {\em SPE Reservoir Evaluation and Engineering}, 4(04):308--317,
  August 2001.

\bibitem{Hoang:1213}
T.~Hoang, J.~Jaffr{\'e}, C.~Japhet, M.~Kern, and J.E. Roberts.
\newblock Space-time domain decomposition methods for diffusion problems in
  mixed formulations.
\newblock {\em SIAM Journal on Numerical Analysis}, 51(6):3532--3559, December
  2013.

\bibitem{Hoang:0717}
T.~Hoang, C.~Japhet, M.~Kern, and J.E. Roberts.
\newblock Space-time domain decomposition for advection--diffusion problems in
  mixed formulations.
\newblock {\em Mathematics and Computers in Simulation}, 137:366--389, July
  2017.

\bibitem{Hughes:0288}
T.J.R. Hughes and G.M. Hulbert.
\newblock Space-time finite element methods for elastodynamics: Formulations
  and error estimates.
\newblock {\em Computer Methods in Applied Mechanics and Engineering},
  66(3):339--363, Feburary 1988.

\bibitem{Hulbert:1290}
G.M. Hulbert and T.J.R. Hughes.
\newblock Space-time finite element methods for second-order hyperbolic
  equations.
\newblock {\em Computer Methods in Applied Mechanics and Engineering},
  84(3):327--348, December 1990.

\bibitem{Kocher:15}
U.~K{\"o}cher.
\newblock {\em Variational Space-Time Methods for the Elastic Wave Equation and
  the Diffusion Equation}.
\newblock PhD thesis, Helmut-Schmidt-University, 2015.

\bibitem{Kocher:1114}
U.~K{\"o}cher and M.~Bause.
\newblock Variational space--time methods for the wave equation.
\newblock {\em Journal of Scientific Computing}, 61(2):424--453, November 2014.

\bibitem{Peszy:0306}
M.~Peszy{\'n}ska, M.F. Wheeler, and I.~Yotov.
\newblock Mortar upscaling for multiphase flow in porous media.
\newblock {\em Computational Geosciences}, 6(1):73--100, March 2006.

\bibitem{Singh:1118}
G.~Singh, W.~Leung, and M.F. Wheeler.
\newblock Multiscale methods for model order reduction of non-linear multiphase
  flow problems.
\newblock {\em Computational Geosciences}, pages 1--19, November 2018.

\bibitem{Singh:0918}
G.~Singh and M.F. Wheeler.
\newblock A domain decomposition approach for local mesh refinement in space
  and time.
\newblock In {\em SPE Annual Technical Conference and Exhibition}. Society of
  Petroleum Engineers, September 2018.

\bibitem{Singh:0818}
G.~Singh and M.F. Wheeler.
\newblock A space-time domain decomposition approach using enhanced velocity
  mixed finite element method.
\newblock {\em Journal of Computational Physics}, 374:893--911, December 2018.

\end{thebibliography}

\end{document}